\documentclass[12pt]{article}
\usepackage{amsmath,amsfonts,amssymb,amsthm}
\input xy
\xyoption{all} \makeatother
\newtheorem{thm}{Theorem}[section]
\newtheorem{pro}[thm]{Proposition}
\newtheorem{lem}[thm]{Lemma}

\theoremstyle{definition}
\newtheorem{dfn}[thm]{Definition}
\newtheorem{rem}{Remark}
\numberwithin{equation}{section}
\newcommand{\CC}{\mathcal{C}}        

\newcommand{\CI}{\mathcal{I}}
\newcommand{\CR}{\mathcal{R}}
\newcommand{\cop}{\Delta}           
\newcommand{\ts}{\otimes}           
\newcommand{\sq}{\unskip\nobreak\kern5pt\nobreak\vrule
     height4pt width4pt depth0pt}   
\newcommand{\Z}{\mathbb{Z}}         
 
\title{On matrix type corings, algebra coverings \\
and \v{C}ech cohomology}
\author{Andrzej Sitarz
\thanks{The author acknowledges the Alexander von Humboldt Fellowship
at the Mathematisches Institut der Heinrich-Heine-Universit\"at,
Universit\"atsstrasse 1, 40225 D\"usseldorf, Germany}
\thanks{Partially supported by Polish Government grants
115/E-343/SPB/6.PR UE/DIE 50/2005--2008 and 189/6.PRUE/2007/7 } \\
Institute of Physics, Jagiellonian University, \\
Reymonta 4, 30-059 Krak\'ow, Poland}
\begin{document}
\maketitle
\begin{abstract}
\noindent
We investigate the a matrix-type coring associated
to a complete covering of an algebra, its Amitsur
complex and propose a definition for the related
\v{C}ech cohomology relative to the covering.
\end{abstract}


\section{Introduction}

The definition and examples of algebra coverings, as given by
Calow and Matthes \cite{CaMa} provided a nice set of examples
of canonical corings \cite{BrWr}. In this note, we elaborate
further on the construction and define a version of \v{C}ech
cohomology respective to the covering using the cochain
complexes based on Amitsur complex.

\subsection{Preliminaries}

We recall here the notation and definition. Let $A$ be a unital
algebra and $I_i \subset A$, $i =1,2,\ldots N$, be two-sided
ideals in $A$. We call $A_j = A/I_j$ and $A_{ij} = A/(I_i + I_j)$,
further $B= \bigoplus_i A_i$ and $B' = \bigoplus_{i,j} A_{ij}$.

\begin{dfn}[Definition 2,\cite{CaMa}]
The ideals $\{I_i\}_{i \in \CI}$ are called a covering of
$A$ if their intersection is $\{ 0 \}$. We call the covering
{\em complete} if the following sequence of $A$-bimodules
is exact:
$$
\xymatrix{0 \ar[r] & A \ar[r]^\pi &B \ar[r]^\tau & B'}
$$
where
$$\pi = \bigoplus_i \pi_i, \,\, \pi_i: A \to A_i,$$
and
$$\tau: \bigoplus_i A_i \to \bigoplus_{i,j} A_{ij}:
\tau(\bigoplus_i x_i) = \oplus_{i,j}
\left( \pi^i_{ij} (x_i) - \pi^j_{ij} (x_j), \right),$$
where $\pi^i_{ij}: A/I_i \to A/(I_i+I_j)$.
\end{dfn}

\begin{dfn}
The coring $\CC(A,\{I\}_\CI)$ associated to the covering
$A, \{I_i\}_{i \in \CI}$ is the Sweedler coring of the
algebra extension $A \hookrightarrow B$.
\end{dfn}

\begin{rem}
The coring $\CC$ is a generalization of a matrix coring,
where the $(i,j)$-matrix entry takes values in
$A_i \ts_A A_j$. The latter can be, for unital algebras
identified (as $A$-bimodule) with $A_{ij}$.

If we denote $e_{ij} = \pi^i(1_A) \ts_A \pi^j(1_A)$
then for the coproduct and counit in $\CC$ we have:

$$
\begin{aligned}
\cop e_{ij} &= \pi^i(1_A) \ts_B
\left( \sum_k \pi^k(1_A) \right) \ts_A \pi_j(1_A) \\
&= \sum_k \left( \pi^i(1_A) \ts_A \pi^k(1_A) \right)
\ts_B
\left( \pi^k(1_A) \ts_A \pi^j(1_A) \right) \\
&= \sum_k e_{ik} \ts_B e_{kj}. \\
\varepsilon(e_{ij}) &= \delta_{ij} \pi_i(1_A).
\end{aligned}
$$
\end{rem}

Note that although formally the coproduct formula is the
same as in the case of a matrix coring, only in the case
$A_i \equiv A, \, i \in I$, we recover the usual
definition of the matrix coring.

We have:
\begin{thm}[\cite{BrWr}, 3.4]
The coring $\CC$ is Galois if the covering is complete.
\end{thm}

\section{The \v{C}ech cohomology of noncommutative coverings}

In the classical situation of coverings of a topological
space with open sets, with some additional data, we can
easily introduce the notion of \v{C}ech cochain complex
relative to the covering.

We show that in the noncommutative situation, this is
also possible and we relate the construction to the
Amitsur complex of the coring related to the covering.
However, for the purpose of the construction we need to
adapt the notion of ringed space to the algebra case and
define the corresponding complex.

\subsection{\v{C}ech complex of rings}

Let $\CI$ be a finite set of indices and $P(\CI)$ be the set of
its ordered subsets. We say that for any two ordered subsets
$\zeta,\vartheta$ of $\CI$, $\zeta \subset \vartheta$ if $\zeta$
is a subset of $\vartheta$ and the order in $\zeta$ is inherited
from the order in $\vartheta$. We shall use $|\zeta|$ to denote
cardinality of the subset $\zeta$.

We assume that $\CR$ is a functor from $P(\CI)$ to the category
of unital rings. That is, for any $\vartheta \in P(\CI)$ there
is a unital ring $\CR(\vartheta)$, and for any $\zeta, \vartheta
\in P(\CI)$ such that $\zeta \subset \vartheta$ there is a ring
homomorphism:
$$ r_\vartheta^\zeta: R(\zeta) \to R(\vartheta). $$

We assume that $r_\zeta^\zeta$ is always an identity morphism,
and that for any three $\zeta \subset \eta \subset
\vartheta \subset  P(\CI)$ we have:
\begin{equation}
r^{\eta}_\vartheta \circ r_\zeta^\eta = r_{\zeta}^\eta.
\label{r1}
\end{equation}
We shall identify $R = \CR(\emptyset)$.


Since $\CR(\emptyset) \hookrightarrow \bigoplus_i \CR(\{i\})$,
is an embedding, we have a similar situation of a
ring extension as before and we can associate a
Sweedler coring with the construction, which we shall
call $\CC(\CI)$.

\subsection{The Amitsur complex of a coring}

We recall here the definition of the Amitsur complex for the
canonical Sweedler coring $\CC = B \ts_A B$ \cite{BrWi}. We
use similar notation as before, here $A_\CR = \CR(\emptyset)$
and $B_\CR = \bigoplus_i \CR(\{i\})$.

\begin{dfn}
Let $\{C^n\}_{n=0,\ldots}$ be the following complex:
$$ C^0 = B_\CI, \;\;\;\; C^n = B_\CI^{\otimes_{A_\CR} (n-1)}
\cong \CC^{\otimes_{B_\CR} n},\; n \geq 1,$$
with the map:
$$
\begin{aligned}
d (x) &= 1 \ts_{A_\CR} x - x \ts_{A_\CR} 1, \;\; x \in C^0, \\
d (\omega) &= 1 \ts_{A_\CR} \omega - b'(\omega) +
(-1)^{n+1} \omega \ts_{A_\CR} 1, \;\; \omega \in C^n,\; n\geq 1,
\end{aligned}
$$
where $b'$ is the coboundary of the Hochschild bar complex
constructed with the coring coproduct:
$$ b' (c_0 \ts_{B_\CR} c_1 \ts_{B_\CR} \cdots \ts_{B_\CR} c_n) =
\sum_{i=0}^n (-1)^n c_0 \ts_{B_\CR} c_1 \ts_{B_\CR} \cdots
\Delta c_i \cdots \ts_{B_\CR} c_n.$$

which, for the Sweedler corings is expressed using the
equivalent presentation as:
$$
\begin{aligned}
b'( x_0 \ts_{A_\CR} x_1 \ts_{A_\CR} \cdots \ts_{A_\CR} x_{n+1}) = \\
\sum_{i=0}^n (-1)^n x_0 \ts_{A_\CR} x_1 \ts_{A_\CR} \cdots x_i
\ts_{A_\CR} 1 \ts_{A_\CR} x_{i+1} \cdots \ts_{A_\CR} x_{n+1}.
\end{aligned}
$$
\end{dfn}

We have (see \cite{BrWi} 29.5)
\begin{pro}
The Amitsur complex of a Galois coring is acyclic. Hence the
Amitsur complex of a coring related to a complete covering is
acyclic.
\end{pro}

Therefore, to obtain more information, which arises from the
covering we need to construct a different cochain complex,
associated to the functor $\CR$.

\begin{dfn}
Let $\phi$ be a map between the complex of the coring of
covering $\CC(\CR)$ and the $\Z$-graded module
$\bigoplus_{n \in \Z} S^n$, where
$$ S^n = \bigoplus_{\zeta \in P(\CI), |\zeta|=n} \CR(\zeta). $$
defined as:

$$ \phi(y_{i_1} \ts_{B_\CR} y_{i_2} \cdots \ts_{B_\CR} y_{i_n})
= r^{i_1}_\zeta( \Phi_{I_1}(y_{i_1})) r^{i_2}_\zeta(\Phi_{I_2}(y_{i_2}))
\cdots r^{i_n}_\zeta(\Phi_{I_n}(y_{i_n})), $$
for any $\zeta = \{ i_1, i_2, \ldots, i_n \} \in P(I)$,
$y_{i_k} \in A_{i_k}$, and:
$$ \phi(y_{i_1} \ts_{B_\CR} y_{i_2} \cdots \ts_{B_\CR} y_{i_n}) \equiv 0,$$
for any $\{ i_1, i_2, \ldots, i_n \}$ which are not in
$P(\CI)$\footnote{This happens when at least two of 
the indices are identical.}. The maps $\Phi_{i_k}$ denote 
some chosen ring homomorphisms between $A_{i_k}$ and $\CR(\{i_k\})$,
in fact, we have $\phi(y_i) \equiv \Phi_{i}(y_i)$.
\end{dfn}

\begin{pro}
$S^n$ is a cochain complex, with the coboundary:
map $d': S^n \to S^{n+1}$:
$$ d' x_\zeta = \sum_{\eta, \zeta \subset \eta,
|\eta| = |\zeta|+1} (-1)^{|\eta/\zeta|}
r_\eta^\zeta(x_\zeta),$$
where $|\eta/\zeta|$ denotes position on which the difference
between $\eta$ and $\zeta$ occurs.
Then $\phi$ is a morphism of cochain complexes.
\end{pro}

We call the cochain complex $S(\CR) = \{(S^n, d')\}_n$
the \v{C}ech complex associated to functor $\CR$.

\begin{proof}
First of all, let us see that the map $\phi$ is well-defined
(the definition of $\phi$ is set for tensor products and we
need to check that it remains correct for tensor products
over $B_\CR$.) We have:
$$
\begin{aligned}
\phi( &\cdots y^{i_k} \pi_{i_k}(y) \ts_{B_\CR} y^{i_{k+1}} \cdots ) =
\cdots r^{i_k}_\zeta \left( \Phi_{i_k}( y_{i_k} \pi_{i_k}(y) ) \right)
r^{i_{k+1}}_\zeta \left( \Phi_{i_{k+1}}( y_{i_{k+1}})\right) \cdots \\
&= \cdots r^{i_k}_\zeta
\left( \Phi_{i_k}( y_{i_k}) \Phi_{i_k}(\pi_{i_k}(y)) \right)
r^{i_{k+1}}_\zeta \left( \Phi_{i_{k+1}}( y_{i_{k+1}})\right) \cdots \\
&= \cdots r^{i_k}_\zeta \left( \Phi_{i_k}( y_{i_k}) \right)
r^{i_k}_\zeta \left( \Phi_{i_k}(\pi_{i_k}(y)) \right)
r^{i_{k+1}}_\zeta \left( \Phi_{i_{k+1}}( y_{i_{k+1}})\right) \cdots \\
&= \cdots r^{i_k}_\zeta \left( \Phi_{i_k}( y_{i_k}) \right)
r^{i_{k+1}}_\zeta \left( \Phi_{i_{k+1}}(\pi_{i_{k+1}}(y)) \right)
r^{i_{k+1}}_\zeta \left( \Phi_{i_{k+1}}( y_{i_{k+1}})\right) \cdots \\
&= \cdots r^{i_k}_\zeta \left( \Phi_{i_k}( y_{i_k}) \right)
r^{i_{k+1}}_\zeta \left(
\Phi_{i_{k+1}}(\pi_{i_{k+1}}(y) y_{i_{k+1}}) \right) \cdots \\
&= \phi(\cdots y^{i_k} \ts_{B_\CR} r_{i_{k+1}}(y) y^{i_{k+1}} \cdots ).
\end{aligned}
$$

As the next step, we check that the map $d'$ is a coboundary:
$$
\begin{aligned}
(d')^2 x_\zeta &= d' \left( \sum_{\eta, \zeta \subset \eta,
|\eta| = |\zeta|+1} (-1)^{|\eta/\zeta|}
r_\eta^\zeta(x_\zeta) \right) \\
&= \sum_{\vartheta, \vartheta \subset \zeta, |\zeta| = |\vartheta|+1}
\sum_{\eta, \zeta \subset \eta, |\eta| = |\zeta|+1}
(-1)^{|\vartheta/\zeta|} (-1)^{|\eta/\zeta|}
r^\eta_\vartheta \circ r^\zeta_\eta (x_\zeta) \\
&= \sum_{\vartheta, \vartheta \subset \zeta, |\zeta| = |\vartheta|+2}
\left( (-1)^{|\eta_1/\zeta|+|\vartheta/\eta_1|} +
(-1)^{|\eta_2/\zeta|+|\vartheta/\eta_2|} \right)
r^\zeta_\vartheta (x_\zeta) \\
&=0,
\end{aligned}
$$
where we have used (\ref{r1}) and the fact that for each
$\vartheta \subset \zeta$ such that $|\vartheta|+2 = |\zeta|$
there are two possibilities for $\vartheta \subset \eta_{1,2} \subset
\zeta$ with $|\eta|=|\vartheta|+1$, which we denoted
$\eta_1, \eta_2$. It is easy to verify that:
$$ (-1)^{|\eta_1/\zeta|+|\vartheta/\eta_1|} +
(-1)^{|\eta_2/\zeta|+|\vartheta/\eta_2|} = 0,$$
which ends the verification.

Finally, to verify that $\phi$ is a cochain complex morphism we
calculate now using $x_i = \phi(y_i)$. First:
$$
\begin{aligned}
d'\circ &\phi (x_{i_0} \ts_R x_{i_1} \cdots \ts_R x_{i_n}) =
d' \left( r^{i_0}_\zeta(x_{i_0}) r^{i_1}_\zeta(x_{i_1})
\cdots r^{i_n}_\zeta(x_{i_n}) \right) \\
&= \sum_{i \in I/\zeta} \sum_{k=0}^{n+1} (-1)^k
r^{\eta(k)}_\zeta \left( r^{i_0}_\zeta(x_{i_1}) r^{i_1}_\zeta(x_{i_1})
\cdots r^{i_n}_\zeta(x_{i_n}) \right) \\
&= \sum_{i \in I/\zeta} \sum_{k=0}^{n+1} (-1)^k
r^{i_0}_{\eta(k)}(x_{i_0}) r^{i_1}_{\eta(k)}(x_{i_1}) \cdots r^{i_n}_{\eta(k)}(x_{i_n})
\end{aligned}
$$
where ${\eta(k)}$ denotes the ordered subset of $P(I)$,
containing $\zeta$ and $i \notin \zeta$ on $k\!+\!1$-th place.

On the other hand, calculating $ \phi \circ d$,
$$
\begin{aligned}
\phi \circ & d (x_{i_0} \ts_R x_{i_2} \cdots \ts_R x_{i_n}) = \\
&= \phi \left( \sum_{k=0}^{n+1} (-1)^k \phi(x_{i_0} \ts_R x_{i_1}
\cdots x_{i_{k-1}} \ts_R (\sum_{i=1}^N r_i(1)) \ts_R x_{i_k}
\cdots \ts_R x_{i_n}) \right) \\
&= \sum_{k=0}^{n+1} r^{i_0}_{\eta(k)}(x_{i_0})
r^{i_1}_{\eta(k)}(x_{i_1}) \cdots \left(
\sum_{i \in I/\zeta} r^i_{\eta(k)}(1) \right) \cdots
r^{i_n}_{\eta(k)}(x_{i_n}) \\
&= \sum_{i \in I/\zeta} \sum_{k=0}^{n+1} r^{i_0}_{\eta(k)}(x_{i_0})
r^{i_1}_{\eta(k)}(x_{i_1}) \cdots r^{i_n}_{\eta(k)}(x_{i_n})
\end{aligned}
$$
which ends the proof.
\end{proof}
\begin{rem}
The trivial example of the functor $\CR$ is given by a constant
functor associating ring $R$ to every $\zeta \in P(\CI)$. The
coring $\CC(\CR)$ is then a full-matrix coring.

The homology of the associated
\v{Cech} cochain complex is:
$$ H^0(S) = R, \;\; H^1(S)= 0. $$
\end{rem}

\begin{rem}
If $X$ is a ringed space with the functor $\Psi$ from open
subsets of $X$ to the category of unital rings, and
$U_i$, $i \in \CI$ is a finite covering of $X$, then we
construct $\CR$ by setting:
$$\CR(\{i_1,i_2,\ldots, i_k \}) = \Psi(U_{i_1} \cap U_{i_2}
\cdots \cap U_{i_k}).$$

In the next section we shall discuss the relation of the
constructed complex with the standard \v{C}ech cohomology.
\end{rem}

\subsection{\v{C}ech cohomology of ringed algebras}

One of the examples of the functor $\CR$ for a finite set
$\CI$ was based on the structure of a ringed space and a
finite covering by open sets. We adapt the definition to
the algebraic case:

\begin{dfn}
We call the algebra $A$ {\em ringed} if there exists a functor
$\Phi$ which associates to any ideal $J \subset A$ a ring
$\Phi(J)$, and a ring morphism $\Phi_J: A/J \to R(J)$, such
that if $J_1 \subset J_2$ then the following diagram is
commutative:

$$
\xymatrix{
A/J_1 \ar[d]^\Phi_{J_1} \ar[r]^\pi & A/J_2 \ar[d]^\Phi_{J_2} \\
\phi(J_1) \ar[r]^{\Phi(\pi)} & \Phi(J_2) }
$$
\end{dfn}

\begin{rem}
An example of the ringed algebra is given by taking
$\Phi(J) = A/J$ and $\Phi_J$ the identity map.
\end{rem}

\begin{pro}
Assume now that $A$ is a ringed algebra with the ringed
structure $\Phi$ and $\{ J_i \}_{i \in \CI}$ is a finite
complete covering of $A$. Then by setting for every ordered
subset of $\CI$, $\zeta=\{i_1,\ldots, i_k\} \in P(\CI)$
$$\CR(\zeta) = \Phi(I_{i_1} + \cdots + I_{i_k}), $$
and for the respective homomorphisms
$r^\eta_\zeta: R(\eta) \to R(\zeta)$ to be images under
$\Phi$ of projection morphisms.
we obtain a functor $\CR(\CI)$.
\end{pro}

\begin{dfn}
If $(A, \Phi)$ is a ringed algebra, and $\{I_i\}_{i \in \CI}$
is its finite covering, we define the \v{C}ech cohomology of
$A$ relative to the covering as the homology of the
\v{C}ech cochain complex of the covering coring.
\end{dfn}

\section{Examples and applications}

Let $X$ be a topological space, and ${\cal O}_X$ be a presheaf
of rings on $X$. We denote the functor associating a ring to
any open set $U$ by $\Phi(U)$, and the ring morphisms
corresponding to set inclusion $i_{U,V}: U \hookrightarrow V$ by
$\Phi_U^V : \Phi(V) \to \Phi(U)$.

Let $U_i$, $i \in I$ be a finite covering of $X$. Then from
the presheaf ${\cal O}_X$ we can construct the data of a ringed
space as in the previous section, the coring, the map $\phi$ and
the respective cochain complexes. We have:

\begin{lem}
The homology of the cochain complex $S(\CR)$ for the coring
given by the data $I, {\cal O}_X$, is the \v{C}ech cohomology
of $X$ relative to the covering $\{U_i\}_I$.
\end{lem}

\begin{proof}
By rewriting the definition of $S^n,d'$ using the sets
$U_i$ and the intersections, we explicitly recover the
definition of \v{C}ech homology.
\end{proof}

\begin{rem}
The construction of \v{C}ech cohomology is possible for
any matrix-type coring, by using the canonical functor
ringed structure $\Phi$ and the covering implicitly defined
by the coring construction.
\end{rem}

\section{Conclusions}

The formulation of {\em noncommutative coverings} in the
language of algebra extensions and associated corings has
enabled us to extend the definition of \v{C}ech cohomology
to this setting.

The \v{C}ech cohomology is mostly of combinatorial nature
and strongly depends on the properties of the underlying
noncommutative covering. It would be interesting to know,
whether for some types of noncommutative coverings one
could establish relation between this cohomology and
some other cohomology theory of the algebra.

\end{document}